# On the path structure of a semimartingale arising from monotone probability theory

Alexander C. R. Belton[1]

*School of Mathematics, Applied Mathematics and Statistics, University College, Cork, Ireland.*

**Abstract.** Let $X$ be the unique normal martingale such that $X_0 = 0$ and

$$\mathrm{d}[X]_t = (1 - t - X_{t-})\,\mathrm{d}X_t + \mathrm{d}t$$

and let $Y_t := X_t + t$ for all $t \geq 0$; the semimartingale $Y$ arises in quantum probability, where it is the monotone-independent analogue of the Poisson process. The trajectories of $Y$ are examined and various probabilistic properties are derived; in particular, the level set $\{t \geq 0\colon Y_t = 1\}$ is shown to be non-empty, compact, perfect and of zero Lebesgue measure. The local times of $Y$ are found to be trivial except for that at level 1; consequently, the jumps of $Y$ are not locally summable.

**Résumé.** Soit $X$ l'unique martingale normale telle que $X_0 = 0$ et

$$\mathrm{d}[X]_t = (1 - t - X_{t-})\,\mathrm{d}X_t + \mathrm{d}t$$

et soit $Y_t := X_t + t$ pour tout $t \geq 0$; la semimartingale $Y$ se manifeste dans la théorie des probabilités quantiques, où c'est analogue du processus de Poisson pour l'indépendance monotone. Les trajectoires de $Y$ sont examinées et diverses propriétés probabilistes sont déduites; en particulier, l'ensemble de niveau $\{t \geq 0\colon Y_t = 1\}$ est montré être non vide, compact, parfait et de mesure de Lebesgue nulle. Les temps locaux de $Y$ sont trouvés être triviaux sauf celui au niveau 1; par conséquent les sauts de $Y$ ne sont pas localements sommables.

## 0. Introduction

The first Azéma martingale, that is, the unique (in law) normal martingale $M$ such that $M_0 = 0$ and

$$\mathrm{d}[M]_t = -M_{t-}\,\mathrm{d}M_t + \mathrm{d}t,$$

---

[1]Present address: Department of Mathematics and Statistics, Fylde College, University of Lancaster, Lancaster, LA1 4YF, United Kingdom, E-mail: a.belton@lancaster.ac.uk

has been the subject of much interest since its appearance in [3], Proposition 118 (see, for example, [4, 13] and [17], Section IV.6); it was the first example to be found of a process without independent increments which possesses the chaotic-representation property. It shall henceforth be referred to as *Azéma's martingale*.

From a quantum-stochastic viewpoint, the process $M$ may be obtained by applying Attal's D transform ([1], Section IV) to the Wiener process. Furthermore, thanks to the factorisation of D provided by vacuum-adapted calculus [5], $M$ appears as a natural object in monotone-independent probability theory; the distribution of $M_t$ (the arcsine law) is a central-limit law which plays a rôle analogous to that played by the Gaussian distribution in the classical framework ([16], Theorem 3.1).

The Poisson distribution also occurs as a limit (the *law of small numbers*): if, for all $n \geq 1$, $(x_{n,m})_{m=1}^n$ is a collection of independent, identically distributed random variables and there exists a constant $\lambda > 0$ such that

$$\lim_{n \to \infty} n\mathbb{E}[x_{n,1}^k] = \lambda \quad \forall k \geq 1,$$

then $x_{n,1} + \cdots + x_{n,n}$ converges in distribution to the Poisson law with mean $\lambda$. (A simple proof of this result is provided in Appendix A.) In the case where $x_{n,1}, \ldots, x_{n,n}$ are Bernoulli random variables taking the values 0 and 1 with mean $\lambda/n$, this is simply the Poisson approximation to the binomial distribution ([8], Example 25.2).

A corresponding theorem holds in the monotone set-up ([16], Theorem 4.1), but now the limit distribution is related to the D transform of the standard Poisson process (with intensity 1 and unit jumps) in the same way as the arcsine law and Azéma's martingale are related above [6]. (This result also holds for free probability: see [20], Theorem 4.) The classical process $Y$ which results is such that $Y_t = X_t + t$ for all $t \geq 0$, where $X$ is the unique normal martingale such that $X_0 = 0$ and

$$d[X]_t = (1 - t - X_{t-}) \, dX_t + dt.$$

This article extends the sample-path analysis of $Y$ (and so $X$) which was begun in [7]. Many similarities are found between $Y$ and Azéma's martingale $M$; for example, they are both determined by a random perfect subset of $\mathbb{R}_+$ and a collection of binary choices, one for each interval in that subset's complement. In Section 1 some results from the theory of martingales are recalled; Section 2 defines the processes $X$ and $Y$ and presents their Markov generators. A random time $G_\infty$ after which $Y$ is deterministic is discussed in Section 3: by Proposition 3.1 and Corollary 3.5, $G_\infty < \infty$ almost surely and, in this case,

$$Y_{t+G_\infty} = -W_{-1}(-\exp(-1-t)) \quad \forall t \geq 0,$$

where $W_{-1}$ is a certain branch of the inverse to the function $z \mapsto ze^z$ (see Notation below). In Section 4 the process $X$ is decomposed into an initial waiting time $S_0$ which is exponentially distributed and an independent normal martingale $Z$ which satisfies the same structure equation as $X$ but has the initial condition $Z_0 = 1$; Lemma 4.2 implies that, for all $t \geq 0$,

$$X_t = \begin{cases} -t & \text{if } t \in [0, S_0[, \\ Z_{t-S_0} - S_0 & \text{if } t \in [S_0, \infty[. \end{cases}$$

Explicit formulae are found for the distribution functions of $G_\infty$ and $J$, a random variable analogous to $G_\infty$ but for $Z$ rather than $X$. In Section 5 it is shown that $(H_t := 1 - (Z_t + t)^{-1})_{t \geq 0}$ is a martingale which is related to Azéma's martingale $M$ by a time change; this gives a simple way to find various properties of the level set $\mathcal{U} := \{t \geq 0 \colon Y_t = 1\}$ in Section 6. Finally, Section 7 presents some results on the local times of $Y$. The appendices contain various supplementary results which are not appropriate for the main text.

*0.1. Conventions*

The underlying probability space is denoted $(\Omega, \mathcal{F}, \mathbb{P})$ and is assumed to contain a filtration $(\mathcal{F}_t)_{t \geq 0}$ which generates the $\sigma$-algebra $\mathcal{F}$. This filtration is supposed to *satisfy the usual conditions*: it is right continuous and the initial $\sigma$-algebra $\mathcal{F}_0$ contains all the $\mathbb{P}$-null sets. Each semimartingale which is considered below has *càdlàg* paths (that is, they are right-continuous with left limits) and two processes $(X_t)_{t\geq 0}$ and $(Y_t)_{t\geq 0}$ are taken equal if they are *indistinguishable*: $\mathbb{P}(X_t = Y_t \text{ for all } t \geq 0) = 1$. Any quadratic variation or stochastic integral has value 0 at time 0.

*0.2. Notation*

The expression $\mathbb{1}_P$ is equal to 1 if the proposition $P$ is true and equal to 0 otherwise; the indicator function of a set $A$ is denoted by $1_A$. The set of natural numbers is denoted by $\mathbb{N} := \{1, 2, 3, \ldots\}$, the set of non-negative rational numbers is denoted by $\mathbb{Q}_+$ and the set of non-negative real numbers is denoted by $\mathbb{R}_+$. The branches of the Lambert $W$ function (that is, the multi-valued inverse to the map $z \mapsto z e^z$) which take (some) real values are denoted by $W_0$ and $W_{-1}$, following the conventions of Corless et al. [10]:

$$W_0(0) = 0, \qquad W_0(x) \in [-1, 0[ \quad \text{and} \quad W_{-1}(x) \in ]-\infty, -1] \quad \forall x \in [-\mathrm{e}^{-1}, 0[.$$

If $\varXi$ is a topological space then $\mathcal{B}(\varXi)$ denotes the Borel $\sigma$-algebra on $\varXi$. The integral of the process $X$ by the semimartingale $R$ will be denoted by $\int X_t \, \mathrm{d}R_t$ or $X \cdot R$, as convenient; the differential notation $X_t \, \mathrm{d}R_t$ will also be employed. The process $X$ stopped at $T$ is denoted by $X^T$, that is, $X_t^T := X_{t \wedge T}$ for all $t \geq 0$, where $x \wedge y$ denotes the minimum of $x$ and $y$. For all $x$, the positive part $x^+ := \max\{x, 0\}$, the maximum of $x$ and 0.

## 1. Normal sigma-martingales and time changes

**Remark 1.1.** Let $A \in \mathcal{F}$ be such that $\mathbb{P}(A) > 0$. If $\mathcal{G}$ is a sub-$\sigma$-algebra of $\mathcal{F}$ such that $A \in \mathcal{G}$ then

$$\widetilde{\mathcal{G}} := \{B \subseteq \varOmega: B \cap A \in \mathcal{G}\}$$

is a $\sigma$-algebra containing $\mathcal{G}$; the map $\mathcal{G} \mapsto \widetilde{\mathcal{G}}$ preserves inclusions and arbitrary intersections. If

$$\widetilde{\mathbb{P}} := \mathbb{P}(\cdot | A) : \widetilde{\mathcal{F}} \to [0, 1]; \qquad B \mapsto \frac{\mathbb{P}(B \cap A)}{\mathbb{P}(A)},$$

then $(\varOmega, \widetilde{\mathcal{F}}, \widetilde{\mathbb{P}})$ is a complete probability space; if $(\mathcal{G})_{t \geq 0}$ is a filtration in $(\varOmega, \mathcal{F}, \mathbb{P})$ satisfying the usual conditions then $(\widetilde{\mathcal{G}}_t)_{t \geq 0}$ is a filtration in $(\varOmega, \widetilde{\mathcal{G}}, \widetilde{\mathbb{P}})$ which satisfies them as well.

If $T$ is a stopping time for the filtration $(\mathcal{G}_t)_{t \geq 0}$ then it is also one for $(\widetilde{\mathcal{G}}_t)_{t \geq 0}$ and, if $B \subseteq \varOmega$,

$$B \in \widetilde{\mathcal{G}_T} \iff B \cap A \in \mathcal{G}_T \iff B \cap A \cap \{T \leq t\} \in \mathcal{G}_t \quad \forall t \geq 0,$$
$$\iff B \cap \{T \leq t\} \in \widetilde{\mathcal{G}}_t \, \forall t \geq 0 \iff B \in (\widetilde{\mathcal{G}})_T,$$

so the notation $\widetilde{\mathcal{G}}_T$ is unambiguous.

**Lemma 1.2.** *If $T$ is a stopping time such that $\mathbb{P}(T < \infty) > 0$ and $M$ is a local martingale then $N : t \mapsto \mathbb{1}_{T < \infty}(M_{t+T} - M_T)$ is a local martingale for the conditional probability measure $\widetilde{\mathbb{P}} := \mathbb{P}(\cdot | T < \infty)$ and the filtration $(\widetilde{\mathcal{F}}_{t+T})_{t \geq 0}$, such that*

$$[N]_t = \mathbb{1}_{T < \infty}([M]_{t+T} - [M]_T) \quad \forall t \geq 0.$$

**Proof.** If $T < \infty$ almost surely and $M$ is uniformly integrable then the first part is immediate, by optional sampling ([18], Theorem II.77.5), and holds in general by localisation and conditioning. The second claim may be verified by realising $[N]$ as a limit of sums in the usual manner (see [17], Theorem II.22, for example). $\square$

**Definition 1.3.** *A martingale $M$ is normal if $t \mapsto (M_t - M_0)^2 - t$ is also a martingale. (If $M_0$ is square integrable then this is equivalent to $t \mapsto M_t^2 - t$ being a martingale, but in general it is a weaker condition.)*

**Definition 1.4.** *A semimartingale $M$ is a sigma-martingale if it can be written as $K \cdot N$, where $N$ is a local martingale and $K$ is a predictable, $N$-integrable process. Equivalently, there exists an increasing sequence $(A_n)_{n \geq 1}$ of predictable sets such that $\bigcup_{n \geq 1} A_n = \mathbb{R}_+ \times \Omega$ and $1_{A_n} \cdot M \in H^1$ for all $n \geq 1$, where $H^1$ denotes the Banach space of martingales $M$ with $\|M\|_{H^1} := \mathbb{E}[[M]_\infty^{1/2}] < \infty$. Every local martingale is a sigma-martingale and if $M$ is a sigma-martingale then so is $H \cdot M$ for any predictable, $M$-integrable process $H$. (The class of sigma-martingales, so named by Delbaen and Schachermayer in [11], was introduced by Chou in [9], where it is denoted $(\Sigma_m)$; the equivalence mentioned above is due to Émery ([12], Proposition 2).)*

**Theorem 1.5 ([14]).** *If $M$ is a semimartingale with $M_0 = 0$ then the following are equivalent:*

(i) *$M$ and $t \mapsto M_t^2 - t$ are sigma-martingales;*
(ii) *$M$ and $t \mapsto [M]_t - t$ are sigma-martingales;*
(iii) *$M$ and $t \mapsto M_t^2 - t$ are martingales;*
(iv) *$M$ and $t \mapsto [M]_t - t$ are martingales.*

**Proof.** Since $M^2 - [M] = 2M_- \cdot M$, the equivalence of (i) and (ii) is immediate; it also follows from this that (iv) implies (iii) ([17], Corollary 3 to Theorem II.27). To complete the proof it suffices to show that (ii) implies (iv).

Suppose (ii) holds and let $(A_n)_{n \geq 1}$ be an increasing sequence of predictable sets such that $\bigcup_{n \geq 1} A_n = \mathbb{R}_+ \times \Omega$ and both $1_{A_n} \cdot M \in H^1$ and $1_{A_n} \cdot N \in H^1$ for all $n \geq 1$, where $N : t \mapsto [M]_t - t$. (Note that if $X \in H^1$ and $B$ is a predictable set then $1_B \cdot X \in H^1$.) Let $T$ be a bounded stopping time; since $1_{A_n} \cdot N$ is a martingale,

$$\mathbb{E}[(1_{A_n} \cdot [M])_T] = \mathbb{E}[(1_{A_n} \cdot N)_T] + \mathbb{E}\left[\int_0^T 1_{A_n}\,ds\right] = \mathbb{E}\left[\int_0^T 1_{A_n}\,ds\right] \tag{1}$$

and therefore $\mathbb{E}[[M]_T] = \mathbb{E}[T] < \infty$, by monotone convergence. It follows that $\mathbb{E}[|N|_T] \leq \mathbb{E}[[M]_T] + \mathbb{E}[T] < \infty$ and $\mathbb{E}[N_T] = \mathbb{E}[[M]_T - T] = 0$, so $N$ is a martingale. (Apply [17], Theorem I.21 to $N$ stopped at $t$ for any $t \geq 0$.) Furthermore, since $(1_{A_n \setminus A_m} \cdot [M])_t \leq (1_{A_m^c} \cdot [M])_t$ for all $m \leq n$ and $t \geq 0$, where $A_m^c := (\mathbb{R}_+ \times \Omega) \setminus A_m$, the sequence $(1_{A_n \cap (]0,t] \times \Omega)} \cdot M)_{n \geq 1}$ is Cauchy in $H^2$, so convergent there; it follows (by [17], Theorem IV.32, say) that $M$ stopped at $t$ is an $H^2$-martingale. □

**Theorem 1.6.** *If $M$ is a normal martingale and $T$ is a stopping time such that $\mathbb{P}(T < \infty) > 0$ then $N : t \mapsto \mathbb{1}_{T < \infty}(M_{t+T} - M_T)$ is a normal martingale (for the measure $\widetilde{\mathbb{P}} := \mathbb{P}(\cdot | T < \infty)$ and the filtration $(\widetilde{\mathcal{F}}_{t+T})_{t \geq 0}$).*

**Proof.** As $M$ and $t \mapsto (M_t - M_0)^2 - t$ are local martingales, so are $N$ and

$$\begin{aligned}Q : t &\mapsto \mathbb{1}_{T < \infty}((M_{t+T} - M_0)^2 - (t+T) - (M_T - M_0)^2 + T) \\ &= \mathbb{1}_{T < \infty}((M_{t+T} - M_T)^2 - t + 2(M_T - M_0)(M_{t+T} - M_T)),\end{aligned}$$

by Lemma 1.2. Hence $t \mapsto (N_t - N_0)^2 - t = Q_t - 2\mathbb{1}_{T < \infty}(M_T - M_0)N_t$ is also a local martingale (as local martingales form a module over the algebra of random variables which are measurable with respect to the initial $\sigma$-algebra) and the conclusion follows from Theorem 1.5. □

**Lemma 1.7.** *If $A$ is a right-continuous, increasing process such that $A_0 \geq 0$ and each $A_t$ is a stopping time then $(\mathcal{F}_{A_t})_{t \geq 0}$ is a filtration which satisfies the usual conditions.*

**Proof.** This is a straightforward exercise. □

**Lemma 1.8.** *Let $K$ and $L$ be independent martingales and let $A$ be a continuous, increasing, $(\mathcal{F}_t^K)_{t \geq 0}$-adapted process with $A_0 = 0$ and $A_\infty = \infty$, where $(\mathcal{F}_t^K)_{t \geq 0}$ denotes the smallest filtration satisfying the usual hypotheses to which $K$ is adapted.*

*If $\mathcal{G}_t := \mathcal{F}_\infty^K \vee \mathcal{F}_t^L$ for all $t \geq 0$ then each $A_t$ is a $(\mathcal{G}_t)_{t \geq 0}$-stopping time, $(\mathcal{G}_{A_t})_{t \geq 0}$ is a filtration satisfying the usual conditions, $L_A$ is a $(\mathcal{G}_{A_t})_{t \geq 0}$-local martingale and $[L_A] = [L]_A$. If $H$ is an $(\mathcal{F}_t^L)_{t \geq 0}$-predictable process which is $L$ integrable then $H_A$ is $(\mathcal{G}_{A_t})_{t \geq 0}$ predictable and $L_A$ integrable, with $(H \cdot L)_A = H_A \cdot L_A$.*

*If $\mathcal{H}_t := \mathcal{F}_t^K \vee \mathcal{F}_t^{L_A}$ for all $t \geq 0$ then $\mathcal{H}_t \subseteq \mathcal{G}_{A_t}$ for all $t \geq 0$. If there exist disjoint, $(\mathcal{H}_t)_{t \geq 0}$-predictable sets $B$ and $C$ such that $1_B \cdot [K] = [K]$ and $1_C \cdot [L]_A = [L]_A$ and if $([K] + [L]_A)^{1/2}$ is $(\mathcal{H}_t)_{t \geq 0}$-locally integrable then $K + L_A$ is a $(\mathcal{H}_t)_{t \geq 0}$-local martingale and $[K + L_A] = [K] + [L]_A$.*

**Proof.** This is immediate from Lemmes 1–3 and Théorème 1 of [21]. □

## 2. The processes $X$ and $Y$

**Definition 2.1.** *Let $X$ be the normal martingale which satisfies the (time-inhomogeneous) structure equation*

$$\mathrm{d}[X]_t = (1 - t - X_{t-}) \, \mathrm{d}X_t + \mathrm{d}t$$

*with initial condition $X_0 = 0$ and let $Y_t := X_t + t$ for all $t \geq 0$. (The process $X$ was introduced in [7], where it was constructed from the quantum stochastic analogue of the Poisson process for monotone independence. Existence also follows directly from [23], Théorème 4.0.2; uniqueness (in law) and the chaotic-representation property hold by [2], Corollary 26.) Then $Y_0 = 0$ and*

$$\mathrm{d}[Y]_t = (1 - Y_{t-}) \, \mathrm{d}Y_t + Y_{t-} \, \mathrm{d}t, \tag{2}$$

*which implies that $\Delta Y_t \in \{0, 1 - Y_{t-}\}$ for all $t > 0$. If*

$$G_t := \sup\{s \in [0, t] \colon Y_s = 1\} \in \{-\infty\} \cup \, ]0, t] \tag{3}$$

*then (by [7], Theorem 24)*

$$Y_t = -W_\bullet(-\exp(-1 - t + G_t)) \tag{4}$$

*for all $t \geq 0$, where $W_\bullet = W_{-1}$ if $Y_t \geq 1$ and $W_\bullet = W_0$ if $Y_t \leq 1$; a little more will be said in Proposition 6.3. (It follows from this description of the trajectories that $X$ and $Y$ are uniformly bounded on $[0, t]$ for all $t \geq 0$.)*

**Definition 2.2.** *Let*

$$a \colon \mathbb{R}_+ \to \, ]0, 1]; \qquad t \mapsto -W_0(-\mathrm{e}^{-1-t}),$$
$$b \colon \mathbb{R}_+ \to [1, \infty[; \qquad t \mapsto -W_{-1}(-\mathrm{e}^{-1-t})$$

*and*

$$c \colon \, ]0, \infty[ \, \to \mathbb{R}_+; \qquad t \mapsto b'(t) - a'(t) = \frac{b(t)}{b(t) - 1} + \frac{a(t)}{1 - a(t)}.$$

*Note that $a(0) = b(0) = 1$, both $a$ and $b$ are homeomorphisms (which may be verified by inspecting their derivatives on $]0, \infty[$) and $c(t) \searrow 1$ as $t \to \infty$.*

**Lemma 2.3.** *For all $t \geq 0$ the random variable $Y_t$ is distributed with an atom at $0$ (of mass $e^{-t}$) and a continuous part with support $[a(t), b(t)]$:*

$$\mathbb{P}(Y_t \in A) = \mathbb{1}_{0 \in A} e^{-t} + \frac{1}{\pi} \int_{A \cap [a(t), b(t)]} \operatorname{Im} \frac{1}{W_{-1}(-y e^{t-y})} \, dy \quad \forall A \in \mathcal{B}(\mathbb{R}).$$

**Proof.** See [7], Corollary 17. □

**Remark 2.4.** *The (classical) Poisson process is simpler when uncompensated; similarly, it is easier to work with $Y$ than with $X$. These processes are strongly Markov (by [2], Theorem 37, for example) and Émery's Itô formula ([13], Proposition 2) implies that, if $f : \mathbb{R} \to \mathbb{R}$ is twice continuously differentiable,*

$$f(X_t) = f(0) + \int_0^t g(X_{s-}, s) \, dX_s + \int_0^t h(X_{s-}, s) \, ds \tag{5}$$

*and*

$$f(Y_t) = f(0) + \int_0^t g(Y_{s-}, 0) \, dX_s + \int_0^t (h(Y_{s-}, 0) + f'(Y_{s-})) \, ds \tag{6}$$

*for all $t \geq 0$, where $g, h : \mathbb{R}^2 \to \mathbb{R}$ are such that*

$$g(x, t) = \mathbb{1}_{x \neq 1-t} \frac{f(1-t) - f(x)}{1 - x - t} + \mathbb{1}_{x = 1-t} f'(1-t)$$

*and*

$$h(x, t) = \mathbb{1}_{x \neq 1-t} \frac{f(1-t) - f(x) - (1 - x - t) f'(x)}{(1 - x - t)^2} + \mathbb{1}_{x = 1-t} \frac{1}{2} f''(1-t)$$

*for all $x, t \in \mathbb{R}$. It follows that*

$$\lim_{\varepsilon \to 0+} \frac{1}{\varepsilon} \mathbb{E}[f(X_{t+\varepsilon}) - f(X_t) | \mathcal{F}_t] = (\Gamma_t^X f)(X_t)$$

*and*

$$\lim_{\varepsilon \to 0+} \frac{1}{\varepsilon} \mathbb{E}[f(Y_{t+\varepsilon}) - f(Y_t) | \mathcal{F}_t] = (\Gamma^Y f)(Y_t),$$

*for almost all $t \geq 0$, where*

$$(\Gamma_t^X f)(x) := \begin{cases} \frac{f(1-t) - f(x) - (1 - x - t) f'(x)}{(1 - x - t)^2} & \text{if } x \neq 1 - t, \\ \frac{1}{2} f''(1-t) & \text{if } x = 1 - t, \end{cases}$$

$$= \mathbb{1}_{x = 1-t} \frac{1}{2} f''(x) + \int_{\mathbb{R} \setminus \{x\}} (f(y) - f(x) - (y - x) f'(x)) \frac{\delta_{1-t}(dy)}{(y - x)^2}, \tag{7}$$

$$(\Gamma^Y f)(x) := \begin{cases} \frac{f(1) - f(x) - x(1 - x) f'(x)}{(1 - x)^2} & \text{if } x \neq 1, \\ \frac{1}{2} f''(1) + f'(1) & \text{if } x = 1, \end{cases}$$

$$= \mathbb{1}_{x = 1} \frac{1}{2} f''(x) + f'(x) + \int_{\mathbb{R} \setminus \{x\}} (f(y) - f(x) - (y - x) f'(x)) \frac{\delta_1(dy)}{(y - x)^2}, \tag{8}$$

*and $\delta_z$ denotes the Dirac measure on $\mathbb{R}$ with support $\{z\}$.*

## 3. The final jump time

**Proposition 3.1.** *If $G_\infty := \sup\{G_t\colon t \geq 0\}$, where $G_t$ is defined in (3), then the random variable $G_\infty$ (the final jump time of $Y$) is almost surely finite and has density*

$$g_\infty : \mathbb{R} \to \mathbb{R}_+; \qquad x \mapsto \mathbb{1}_{x \geq 0} \frac{1}{\pi} \operatorname{Im} \frac{1}{W_{-1}(-\mathrm{e}^{-1+x})}. \tag{9}$$

**Proof.** Note first that $G_t = 1 + t - Y_t + \log Y_t$ for all $t \geq 0$, by (4), so $G_t$ is $\mathcal{F}_t$ measurable. As $t \mapsto G_t$ is increasing, it is elementary to verify that

$$G_\infty = \sup\{G_t\colon t \geq 0\} = \sup\{G_n\colon n \geq 1\} = \lim_{n \to \infty} G_n;$$

in particular, $G_\infty$ is $\mathcal{F}$ measurable. If $t > 0$ then $\mathbb{1}_{G_n \in ]0,t]} \to \mathbb{1}_{G_\infty \in ]0,t]}$, because $G_n \nearrow G_\infty$, and the dominated-convergence theorem implies that

$$\mathbb{P}(G_\infty \in ]0,t]) = \mathbb{E}[\mathbb{1}_{G_\infty \in ]0,t]}] = \lim_{n \to \infty} \mathbb{E}[\mathbb{1}_{G_n \in ]0,t]}] = \lim_{n \to \infty} \mathbb{P}(G_n \in ]0,t]).$$

Since $\mathbb{P}(G_\infty = -\infty) = \mathbb{P}(Y \equiv 0) \leq \mathbb{P}(Y_t = 0) = \mathrm{e}^{-t} \to 0$ as $t \to \infty$, it follows that $\mathbb{P}(G_\infty = -\infty) = 0$ and

$$\mathbb{P}(G_\infty \leq t) = \lim_{n \to \infty} \mathbb{P}(G_n \in ]0,t])$$

for all $t \geq 0$. If $n \geq 1$ and $t \in [0,n]$ then

$$0 < 1 + n - Y_n + \log Y_n \leq t \quad \iff \quad -\mathrm{e}^{-1-n} > -Y_n \exp(-Y_n) \geq -\mathrm{e}^{-1-n+t}$$
$$\iff \quad Y_n \in \,]a(n), a(n-t)] \cup [b(n-t), b(n)[$$

and, by Lemma 2.3,

$$\gamma_n(t) := \mathbb{P}(G_n \in ]0,t]) = \frac{1}{\pi} \int_{]a(n),a(n-t)] \cup [b(n-t),b(n)[} \operatorname{Im} \frac{1}{W_{-1}(-y\mathrm{e}^{n-y})} \,\mathrm{d}y.$$

Note that $\gamma_n$ is continuously differentiable on $[0,n[$, with

$$\gamma'_n(s) = \frac{1}{\pi} \operatorname{Im} \frac{1}{W_{-1}(-\mathrm{e}^{-1+s})} (b'(n-s) - a'(n-s)) = c(n-s) g_\infty(s)$$

for all $s \in [0,n[$. If $n > t$ and $s \in [0,t]$ then, by the remarks in Definition 2.2, $\gamma'_n(s) \searrow g_\infty(s)$ as $n \to \infty$ and the monotone-convergence theorem implies that

$$\lim_{n \to \infty} \int_0^t \gamma'_n(s) \,\mathrm{d}s = \int_0^t g_\infty(s) \,\mathrm{d}s \quad \forall t \geq 0.$$

This gives the result, because $\int_0^\infty g_\infty(s) \,\mathrm{d}s = 1$ (by Proposition B.1). $\square$

**Remark 3.2.** *It follows from Proposition 3.1 that $\mathbb{E}[G_\infty] = \infty$; a proof is given in Proposition B.1.*

**Remark 3.3.** *Calling $G_\infty$ the final jump time is perhaps a little misleading, since it is not a stopping time; it is, however, almost surely the limit of a sequence of jump times. (See Corollary 6.2 and Corollary 6.4.)*

**Proposition 3.4.** $\lim_{t \to \infty} \mathbb{P}(Y_t \leq 1) = 0.$

**Proof.** By Lemma 2.3,

$$\mathbb{P}(Y_t \leq 1) = e^{-t} + \frac{1}{\pi} \int_{a(t)}^{1} \operatorname{Im} \frac{1}{W_{-1}(-y e^{t-y})} \, dy \quad \forall t \geq 0. \tag{10}$$

If $y \in \,]0,1]$ then there exists $x \in [0, \infty[$ such that $y = a(x)$, and if $t \geq x$ then

$$\operatorname{Im} \frac{1}{W_{-1}(-y e^{t-y})} = \operatorname{Im} \frac{1}{W_{-1}(-e^{-1+t-x})} = \pi g_\infty(t-x) \to 0$$

as $t \to \infty$. (This last claim follows from Proposition B.1.) Furthermore, as $g_\infty$ is bounded, the integrand in (10) is bounded uniformly in $y$ and $t$, so the result follows from the dominated-convergence theorem. $\square$

**Corollary 3.5.** *As $t \to \infty$, the process $Y_t \to \infty$ almost surely.*

**Proof.** If $G_\infty < \infty$ then, as $t \to \infty$, either $Y_t \to 0$ or $Y_t \to \infty$; furthermore,

$$\{G_\infty < \infty\} \cap \left\{ \lim_{t \to \infty} Y_t = \infty \right\} = \{G_\infty < \infty\} \cap \bigcap_{n=1}^{\infty} \bigcup_{m=n}^{\infty} \{Y_m > 1\}.$$

Since $\mathbb{P}(G_\infty < \infty) = 1$ and $\mathbb{P}(Y_n \leq 1) \to 0$ as $n \to \infty$, it follows that

$$\mathbb{P}\left(\lim_{t \to \infty} Y_t = \infty\right) \geq \limsup_{n \to \infty} \mathbb{P}(Y_n > 1) = 1 - \lim_{n \to \infty} \mathbb{P}(Y_n \leq 1) = 1.$$

(The inequality in the previous line holds by [8], Theorem 4.1(i).) $\square$

## 4. The active period

**Proposition 4.1.** *The stopping time $S_0 := \inf\{t > 0 \colon Y_t = 1\}$ is exponentially distributed and has mean 1.*

**Proof.** Note that $Y_t = 0$ only if $Y_s = 0$ for all $s \in [0, t[$, by (4); the claim now follows from Lemma 2.3. $\square$

**Lemma 4.2.** *If $Z_t := X_{t+S_0} + S_0$ for all $t \geq 0$ then $Z$ is a normal martingale for the filtration $(\mathcal{F}_{t+S_0})_{t \geq 0}$ such that $Z_0 = 1$, which satisfies the structure equation*

$$d[Z]_t = dt + (1 - t - Z_{t-}) \, dZ_t \tag{11}$$

*and which is independent of $\mathcal{F}_{S_0}$.*

**Proof.** As $Z_t = X_{t+S_0} - X_{S_0} + 1$ for all $t \geq 0$, Theorem 1.6 implies that $Z$ is a normal martingale. Furthermore,

$$[Z]_t = [X]_{t+S_0} - [X]_{S_0} = \int_{S_0}^{t+S_0} (1 - r - X_{r-}) \, dX_r = \int_0^t (1 - s - Z_{s-}) \, dZ_s$$

for all $t \geq 0$. (The first equality is a consequence of Lemma 1.2; the last may be shown by expressing the integrals as the limit of Riemann sums, as in [17], Theorem II.21, for example.) It now follows from [2], Theorem 25, that, for all $t \geq 0$, the law of $Z_t$ conditional on $\mathcal{F}_{S_0}$ depends only on the initial value $Z_0 = 1$ and the coefficient functions $\alpha \colon s \mapsto 1 - s$ and $\beta \equiv -1$ restricted to $[0, t]$, so $Z_t$ is independent of $\mathcal{F}_{S_0}$. $\square$

*Remark 4.3.* If $t \geq 0$ then

$$Z_t + t = Y_{t+S_0} = -W_\bullet(-\exp(-1 - (t + S_0) + G_{t+S_0})) \in [a(t), b(t)],$$

since $G_{t+S_0} \geq S_0$. Consequently, $Z$ is uniformly bounded on $[0, t]$ for all $t \geq 0$.

**Remark 4.4.** Let $m_n(t) := \mathbb{E}[(Z_t + t)^n]$ for all $n \geq 1$ and $t \geq 0$, where $Z$ is as in Lemma 4.2. It may be shown using Émery's Itô formula ([13], Proposition 2 and the subsequent remark) that

$$m_n(t) - m_{n-1}(t) = n \int_0^t m_{n-1}(s)\,\mathrm{d}s \tag{12}$$

for all $n \geq 1$ and $t \geq 0$ (where $m_0 \equiv 1$). Hence (compare [6], Section 4)

$$\widehat{m}_n(p) = p^{-1} \prod_{j=1}^n (1 + jp^{-1})$$

if $n \geq 1$, where $\widehat{f}$ denotes the Laplace transform of $f$, and so

$$m_n(t) = 1 + \sum_{k=1}^n \left( \sum_{1 \leq j_1 < \cdots < j_k \leq n} j_1 \cdots j_k \right) \frac{t^k}{k!} = \sum_{k=0}^n \begin{bmatrix} n+1 \\ n+1-k \end{bmatrix} \frac{t^k}{k!} \tag{13}$$

for all $t \geq 0$, where $\begin{bmatrix} \cdot \\ \cdot \end{bmatrix}$ denotes the unsigned Stirling numbers of the first kind [15]. (The final identity holds by [7], Proposition 3 and Remark 6, for example.)

**Theorem 4.5.** If $t > 0$ then $Z_t + t = Y_{t+S_0}$ is continuously distributed, with density

$$f_{Z_t+t} : \mathbb{R} \to \mathbb{R}_+; \qquad z \mapsto \mathbb{1}_{z \in [a(t),b(t)]} \frac{1}{\pi} \operatorname{Im} \frac{1}{1 + W_{-1}(-z\mathrm{e}^{t-z})}. \tag{14}$$

**Proof.** Let $x \geq 0$. Since $Y_t = \mathbb{1}_{t \geq S_0}(Z_{(t-S_0)^+} + t - S_0)$ for all $t \geq 0$, it follows that

$$\mathbb{P}(0 < Y_t \leq x) = \mathbb{P}(S_0 \leq t \text{ and } Z_{(t-S_0)^+} + t - S_0 \leq x)$$

$$= \int_0^t \int_{-\infty}^{x-t+s} \mathrm{d}F_{Z_{t-s}}(z) \mathrm{e}^{-s}\,\mathrm{d}s$$

$$= \mathrm{e}^{-t} \int_0^t \int_{-\infty}^{x-u} \mathrm{d}F_{Z_u}(z) \mathrm{e}^u\,\mathrm{d}u,$$

where $F_V$ denotes the distribution function of the random variable $V$. (For the second equality, note that

$$\mathbb{E}[\mathbb{1}_{S_0 \leq t} \mathbb{1}_{Z_{(t-S_0)^+} + t - S_0 \leq x}] = \mathbb{E}[\mathbb{1}_{S_0 \leq t} \mathbb{E}[\mathbb{1}_{Z_{(t-S_0)^+} + t - S_0 \leq x} | \mathcal{F}_{S_0}]]$$

$$= \mathbb{E}[\mathbb{1}_{S_0 \leq t} \mathbb{E}[\mathbb{1}_{Z_{(t-s)^+} + t - s \leq x}]|_{s=S_0}],$$

since $Z$ is independent of $\mathcal{F}_{S_0}$.) Hence

$$\mathbb{P}(Z_t + t \leq x) = \mathrm{e}^{-t} \frac{\mathrm{d}}{\mathrm{d}t}(\mathrm{e}^t \mathbb{P}(0 < Y_t \leq x)) = \mathbb{P}(0 < Y_t \leq x) + \frac{\mathrm{d}}{\mathrm{d}t}\mathbb{P}(0 < Y_t \leq x).$$

Thus if $t > 0$ then either $x \leq a(t)$, so that $\mathbb{P}(Z_t + t \leq x) = 0$, or $x \geq b(t)$, whence $\mathbb{P}(Z_t + t \leq x) = 1 - \mathrm{e}^{-t} + \mathrm{e}^{-t} = 1$, or $x \in {]}a(t), b(t){[}$, in which case

$$\pi \mathbb{P}(Z_t + t \leq x) = \int_{a(t)}^x \operatorname{Im} \frac{1}{W_{-1}(-z\mathrm{e}^{t-z})}\,\mathrm{d}z - a'(t) \operatorname{Im} \frac{1}{W_{-1}(-a(t)\mathrm{e}^{t-a(t)})} + \int_{a(t)}^x \frac{\partial}{\partial t} \operatorname{Im} \frac{1}{W_{-1}(-z\mathrm{e}^{t-z})}\,\mathrm{d}z$$

$$= \int_{a(t)}^x \operatorname{Im} \frac{1}{1 + W_{-1}(-z\mathrm{e}^{t-z})}\,\mathrm{d}z,$$

as claimed. (This formal working is a little awkward to justify: a rigorous proof is provided by Proposition C.2.) □

**Proposition 4.6.** *The random variables $S_0$ and $J := G_\infty - S_0$ are independent and $J$ is continuous, with density*

$$f_J : \mathbb{R} \to \mathbb{R}_+; \qquad x \mapsto \mathbb{1}_{x>0} \frac{1}{\pi} \operatorname{Im} \frac{1}{1 + W_{-1}(-e^{-1+x})}. \tag{15}$$

**Proof.** To see that $S_0$ and $J$ are independent, note first that

$$J = \lim_{n\to\infty} G_{n+S_0} - S_0 = \lim_{n\to\infty} (1 - Z_n + \log(Z_n + n))$$

almost surely, where $(Z_t)_{t \geq 0}$ is as defined in Lemma 4.2, which implies that $G_{n+S_0} - S_0$ is independent of $S_0$ for all $n \geq 1$ and, therefore, so is $J$.

If $F_J(z) := \mathbb{P}(J \leq z)$ for all $z \in \mathbb{R}$ then, by independence and Proposition 4.1,

$$\int_{-\infty}^z g_\infty(w)\, \mathrm{d}w = \mathbb{P}(J + S_0 \leq z) = \iint_{\{(x,y)\in\mathbb{R}^2:\, x+y\leq z\}} \mathrm{d}F_J(x)\, \mathbb{1}_{y\geq 0} e^{-y}\, \mathrm{d}y$$

$$= \int_{-\infty}^z e^{-v} \int_{-\infty}^v e^u\, \mathrm{d}F_J(u)\, \mathrm{d}v$$

for all $z \in \mathbb{R}$, using the substitution $(u, v) = (x, x+y)$. Thus, for almost all $v \in \mathbb{R}$,

$$g_\infty(v) = e^{-v} \int_{-\infty}^v e^u\, \mathrm{d}F_J(v);$$

in fact, this holds for all $v \in \mathbb{R}$, as both functions are continuous, and, since $g_\infty(0) = 0$,

$$g_\infty(t) = e^{-t} \int_0^t e^s\, \mathrm{d}F_J(s) \quad \forall t \geq 0.$$

Now $g_\infty$ is continuously differentiable on $\mathbb{R} \setminus \{0\}$ and $f_J(x) = g_\infty(x) + \mathbb{1}_{x\neq 0} g'_\infty(x)$, so if $0 < \varepsilon < t$ then integration by parts yields the equality

$$\int_\varepsilon^t e^s f_J(s)\, \mathrm{d}s = e^t g_\infty(t) - e^\varepsilon g_\infty(\varepsilon) \to \int_0^t e^s\, \mathrm{d}F_J(s) \quad \text{as } \varepsilon \to 0+.$$

Hence $\int_0^t e^s f_J(s)\, \mathrm{d}s$ exists for all $t \geq 0$ (as does $\int_0^t f_J(s)\, \mathrm{d}s$, by comparison) and

$$\mu : \mathcal{B}(\mathbb{R}_+) \to \mathbb{R}_+; \qquad A \mapsto \int_A e^s\, \mathrm{d}F_J(s) = \int_A e^s f_J(s)\, \mathrm{d}s$$

is a positive Borel measure on $\mathbb{R}_+$; by [19], Theorem 1.29,

$$\int_0^t f_J(s)\, \mathrm{d}s = \int_0^t e^{-s}\, \mathrm{d}\mu(s) = \int_0^t \mathrm{d}F_J(s) = F_J(t) - F_J(0)$$

for all $t \geq 0$ and

$$1 = \lim_{t\to\infty} F_J(t) = F_J(0) + \int_0^\infty g_\infty(s)\, \mathrm{d}s + \lim_{t\to\infty} g_\infty(t) = F_J(0) + 1,$$

by Proposition B.1, so $F_J(0) = 0$. The result follows. $\square$

**Remark 4.7.** *The distribution of $J$ may also be found by imitating the proof of Proposition 3.1, with $Z_t + t$ replacing $Y_t$, since $J$ has the same relationship to $Z$ as $G_\infty$ does to $X$.*

**Proposition 4.8.** *If $t \geq 0$ then*

$$\mathbb{P}(G_\infty \leq t) = -\frac{1}{\pi} \operatorname{Im}\left(W_{-1}(-\mathrm{e}^{-1+t}) + \frac{1}{W_{-1}(-\mathrm{e}^{-1+t})}\right) \tag{16}$$

*and*

$$\mathbb{P}(J \leq t) = -\frac{1}{\pi} \operatorname{Im} W_{-1}(-\mathrm{e}^{-1+t}) = \mathbb{P}(G_\infty \leq t) + g_\infty(t). \tag{17}$$

**Proof.** These follow immediately from the identities

$$\int_0^t \frac{1}{W_{-1}(-\mathrm{e}^{-1+x})} \, \mathrm{d}x = t - \frac{(1 + W_{-1}(-\mathrm{e}^{-1+t}))^2}{W_{-1}(-\mathrm{e}^{-1+t})}$$

and

$$\int_0^t \frac{1}{1 + W_{-1}(-\mathrm{e}^{-1+x})} \, \mathrm{d}x = t - W_{-1}(-\mathrm{e}^{-1+t}) - 1,$$

which are valid for all $t \geq 0$ and may be verified by differentiation. For brevity, let $w = W_{-1}(-\mathrm{e}^{-1+t})$ and $w' = W'_{-1}(-\mathrm{e}^{-1+t})$; note that $\mathrm{d}w/\mathrm{d}t = -\mathrm{e}^{-1+t}w'$ and $-\mathrm{e}^{-1+t}(1+w)w' = w$, whence

$$\frac{\mathrm{d}}{\mathrm{d}t}\left(t - \frac{(1+w)^2}{w}\right) = 1 - \frac{-2\mathrm{e}^{-1+t}(1+w)w'w + \mathrm{e}^{-1+t}w'(1+w)^2}{w^2}$$

$$= 1 - \frac{-\mathrm{e}^{-1+t}(1+w)w'(2w - (1+w))}{w^2} = 1 - \frac{w-1}{w} = \frac{1}{w}$$

and, if $t > 0$,

$$\frac{\mathrm{d}}{\mathrm{d}t}(t - w) = 1 + \mathrm{e}^{-1+t}w' = 1 - \frac{w}{1+w} = \frac{1}{1+w},$$

as required. (To see the existence of $\int_0^t 1/(1 + W_{-1}(-\mathrm{e}^{-1+x})) \, \mathrm{d}x$, note that if $t \geq \varepsilon > 0$ then, letting $W_{-1}(-\mathrm{e}^{-1+x}) = -v \cot v + \mathrm{i}v$, where $v \in {]}-\pi, 0[$,

$$\int_\varepsilon^t \left|\frac{1}{1 + W_{-1}(-\mathrm{e}^{-1+x})}\right| \mathrm{d}x = \int_{v(t)}^{v(\varepsilon)} \sqrt{\frac{1 - 2v\cot v + v^2 \operatorname{cosec}^2 v}{v^2}} \, \mathrm{d}v$$

and the function $v \mapsto (1 - 2v\cot v + v^2 \operatorname{cosec}^2 v)/v^2$ is continuous on ${]}-\pi, 0[$ with limit $1$ as $v \to 0-$.) □

## 5. La martingale cachée

The martingale $H$ discussed in this section was discovered by Émery [14].

**Theorem 5.1.** *If $H_t := 1 - (Z_t + t)^{-1}$ for all $t \geq 0$ then $H$ is a martingale such that $H_0 = 0$,*

$$\mathrm{d}[H]_t = (1 - H_{t-})^2 \, \mathrm{d}t - H_{t-} \, \mathrm{d}H_t \tag{18}$$

*and $H_t \to H_\infty := 1$ almost surely as $t \to \infty$.*

**Proof.** If $t \geq 0$ and $\mathcal{E}(-Z)$ denotes the Doléans-Dade exponential of the normal martingale $-Z$ then $\mathcal{E}(-Z)$ is square integrable on $[0,t]$ for all $t \geq 0$ and (11) implies that

$$(Z_t + t)\mathcal{E}(-Z)_t$$
$$= Z_t + t - \int_0^t (\mathrm{d}Z_s + \mathrm{d}s) \int_0^t \mathcal{E}(-Z)_{s-} \, \mathrm{d}Z_s$$
$$= Z_t + t - \int_0^t (Z_{s-} + s)\mathcal{E}(-Z)_{s-} \, \mathrm{d}Z_s - \int_0^t (1 - \mathcal{E}(-Z)_{s-})(\mathrm{d}Z_s + \mathrm{d}s) - \int_0^t \mathcal{E}(-Z)_{s-} \, \mathrm{d}[Z]_s$$
$$= 1.$$

Thus $H = 1 - \mathcal{E}(-Z)$ is a martingale and $\mathrm{d}H_t = \mathcal{E}(-Z)_{t-} \, \mathrm{d}Z_t = (1 - H_{t-}) \, \mathrm{d}Z_t$, whence

$$\mathrm{d}[H]_t = (1 - H_{t-})^2 \, \mathrm{d}[Z]_t$$
$$= (1 - H_{t-})^2 (\mathrm{d}t + (1 - (1 - H_{t-})^{-1}) \, \mathrm{d}Z_t)$$
$$= (1 - H_{t-})^2 \, \mathrm{d}t - H_{t-} \, \mathrm{d}H_t,$$

as claimed. Since $Y_t \to \infty$ almost surely as $t \to \infty$, by Corollary 3.5, so does $Z_t + t = Y_{t+S_0}$, and the final claim follows. $\square$

**Remark 5.2.** As $H_t = 0$ if and only if $Z_t + t = 1$,

$$\mathcal{U} := \{t \geq 0 \colon Y_t = 1\} = \{s + S_0 \colon Y_{s+S_0} = 1\} = \{s + S_0 \colon H_s = 0\};$$

the structure of $\mathcal{U}$ is determined by the zero set of $H$.

**Definition 5.3.** Let

$$\tau \colon \mathbb{R}_+ \times \Omega \to \mathbb{R}_+; \qquad (t, \omega) \mapsto \tau_t(\omega) := \int_0^t (1 - H_{s-}(\omega))^2 \, \mathrm{d}s$$

and note that $\tau$ is adapted to the filtration $(\mathcal{F}_t)_{t \geq 0}$ and has paths which are continuous, strictly increasing and bi-Lipschitzian on any compact subinterval of $\mathbb{R}_+$, since the derivative

$$\tau'_t = (1 - H_{t-})^2 = (Z_{t-} + t)^{-2} \in [b(t)^{-2}, a(t)^{-2}]$$

for all $t \geq 0$. Let

$$\tau_\infty := \int_0^\infty (1 - H_{s-})^2 \, \mathrm{d}s \in \,]0, \infty]$$

and extend $\tau^{-1}$ (defined pathwise) to all of $\mathbb{R}_+$ by letting $\tau_s^{-1} := \infty$ for all $s \in [\tau_\infty, \infty[$. If $s \geq 0$ then $\{\tau_s^{-1} \leq t\} = \{s \leq \tau_t\} \in \mathcal{F}_t$ for all $t \geq 0$, so $\tau_s^{-1}$ is an $(\mathcal{F}_t)_{t \geq 0}$ stopping time. Thus $(\mathcal{G}_s := \mathcal{F}_{\tau_s^{-1}})_{s \geq 0}$ is a filtration which satisfies the usual conditions, by Lemma 1.7.

**Proposition 5.4.** *The process $K = (K_s := H_{\tau_s^{-1}})_{s \geq 0}$ is a martingale for the filtration $(\mathcal{G}_s)_{s \geq 0}$ and satisfies the equation*

$$[K]_s = s \wedge \tau_\infty - \int_0^s K_{r-} \, \mathrm{d}K_r \quad \forall s \geq 0. \tag{19}$$

**Proof.** Fix $s \geq 0$; as $\tau_s^{-1}$ is an $(\mathcal{F}_t)_{t \geq 0}$ stopping time, $H^{\tau_s^{-1}}$ is a martingale for this filtration ([18], Theorem II.77.4). Let $(T_n)_{n \geq 1}$ be an increasing sequence of stopping times which reduces the local martingale $H_- \cdot H$ and note that

$$\mathbb{E}[\tau_{\tau_s^{-1} \wedge T_n} - [H]_{\tau_s^{-1} \wedge T_n}] = \mathbb{E}[(H_- \cdot H)_{\tau_s^{-1}}^{T_n}] = 0,$$

by the optional-sampling theorem. As $\tau$ is increasing, the monotone-convergence theorem implies that

$$\mathbb{E}[s \wedge \tau_\infty] = \lim_{n \to \infty} \mathbb{E}[\tau_{\tau_s^{-1} \wedge T_n}] = \lim_{n \to \infty} \mathbb{E}[[H]_{\tau_s^{-1} \wedge T_n}] = \mathbb{E}[[H^{\tau_s^{-1}}]_\infty],$$

so $H^{\tau_s^{-1}}$ is a square-integrable martingale ([17], Corollary 4 to Theorem II.27). Hence $K$ is a martingale, by a further application of the optional-sampling theorem: if $0 \leq r \leq s$ then

$$\mathbb{E}[K_s | \mathcal{G}_r] = \mathbb{E}[H_\infty^{\tau_s^{-1}} | \mathcal{F}_{\tau_r^{-1}}] = H_{\tau_r^{-1}} = K_r.$$

Moreover,

$$\int_0^s K_{r-} \, dK_r = \int_0^{\tau_s^{-1}} K_{\tau_r -} \, dH_r = \int_0^{\tau_s^{-1}} H_{r-} \, dH_r$$

(which follows from [17], Theorem II.21, for example), so

$$[K]_s = K_s^2 - K_0^2 - 2 \int_0^s K_{r-} \, dK_r = H_{\tau_s^{-1}}^2 - H_0^2 - 2 \int_0^{\tau_s^{-1}} H_{r-} \, dH_r = [H]_{\tau_s^{-1}}$$

and this equals

$$\tau_{\tau_s^{-1}} - \int_0^{\tau_s^{-1}} H_{r-} \, dH_r = s \wedge \tau_\infty - \int_0^s K_{r-} \, dK_r.$$

□

**Theorem 5.5.** *Let $M$ be Azéma's martingale, that is, the normal martingale such that $M_0 = 0$ and*

$$d[M]_t = dt - M_{t-} \, dM_t.$$

*If $T := \inf\{t \geq 0 : M_t = 1\}$ then $M^T$ and $K$ are identical in law.*

**Proof.** Let $L$ be a normal martingale which is independent of $K$ such that $L_0 = 1$ and

$$d[L]_t = dt - L_{t-} \, dL_t,$$

that is, $L$ is an Azéma's martingale started at 1; existence of such follows from [13], Proposition 5. For all $t \geq 0$, let

$$P_t := \mathbb{1}_{t \in [0, \tau_\infty[} K_t + \mathbb{1}_{t \in [\tau_\infty, \infty[} L_{t - \tau_\infty} = K_t + L_{(t - \tau_\infty)^+} - 1.$$

In the notation of Lemma 1.8, $\tau_\infty = \inf\{t \geq 0 : K_t = 1\}$ is a $(\mathcal{F}_t^K)_{t \geq 0}$-stopping time, so $]0, \tau_\infty]$ is $(\mathcal{F}_t^K)_{t \geq 0}$ predictable and $\mathbb{1}_{]0, \tau_\infty]} \cdot [K] = [K]$ (since $K = K^{\tau_\infty}$) whereas $\mathbb{1}_{]0, \tau_\infty]} \cdot [L_A] = 0$, if $A_t := (t - \tau_\infty)^+$, because $(L_A)_t^{\tau_\infty} = L_{A_{t \wedge \tau_\infty}} = 0$ for all $t \geq 0$. Since $[K]_t = 2(t \wedge \tau_\infty) - K_t^2 \leq 2t$ and $[L]_{A_t} = 2A_t - L_{A_t}^2 \leq 2t$, Lemma 1.8 implies that $P = K + L_A - 1$ is a local martingale such that $P_0 = 0$ and

$$[P]_t = [K]_t + [L]_{A_t} = t - (K_- \cdot K)_t - (L_- \cdot L)_{A_t}.$$

However,

$$[P] = [K] + [L_A]$$
$$= K^2 - 2K_- \cdot K + L_A^2 - 1 - 2L_{A-} \cdot L_A$$
$$= (K + L_A - 1)^2 + 2K + 2L_A - 2 - 2KL_A - 2K_- \cdot K - 2(L_- \cdot L)_A$$

and $KL_A = P$, so

$$P^2 - 2P_- \cdot P = [P] = P^2 - 2K_- \cdot K - 2(L_- \cdot L)_A.$$

Thus $[P]_t = t - (P_- \cdot P)_t$, so $P$ is a normal martingale, by Theorem 1.5, and, by uniqueness ([13], Proposition 6), $P$ is equal to $M$ in law. Since $\tau_\infty = \inf\{t \geq 0 \colon P_t = 1\}$, the processes $K = P^{\tau_\infty}$ and $M^T$ are identical in law, as claimed. □

## 6. The level set $\mathcal{U}$

The level set

$$\mathcal{U} = \{t + S_0 \colon H_t = 0\} = \tau^{-1}(\{s \in [0, \tau_\infty[ \colon K_s = 0\}) + S_0,$$

where $\tau$ is a homeomorphism between $\mathbb{R}_+$ and $[0, \tau_\infty[$ which is bi-Lipschitzian on compact subintervals. This observation, together with Theorem 5.5, leads immediately to the following theorem, thanks to well-known properties of the zero set of Azéma's martingale (or rather, by [17], Section IV.6, properties of the zero set of Brownian motion: see [8], Theorem 37.4 and [24]).

**Theorem 6.1.** *The set $\mathcal{U} := \{t \geq 0 \colon Y_t = 1\}$ is almost surely non-empty, perfect (that is, closed and without isolated points), compact and of zero Lebesgue measure. If $a > 0$ then $\mathcal{U} \cap [S_0, S_0 + a]$ has Hausdorff dimension $1/2$.*

**Corollary 6.2.** *If $T$ is a stopping time then $\mathbb{P}(G_\infty = T) = 0$. In particular, the final jump time $G_\infty$ is not a stopping time.*

**Proof.** If $T$ is a stopping time then so is $T' = \mathbb{1}_{Y_T=1} T + \mathbb{1}_{Y_T \neq 1} \infty$; let $Z'_t := \mathbb{1}_{T' < \infty}(X_{t+T'} - X_{T'} + 1)$ for all $t \geq 0$. Conditional on $T' < \infty$, it holds that $Z'_0 = 1$ and, working as in the proof of Lemma 4.2,

$$\mathrm{d}[Z']_t = (1 - t - Z'_{t-}) \mathrm{d} Z'_t + \mathrm{d}t,$$

so $Z'$ is identical in law to $Z$. In particular, the set $\mathcal{U} \cap ]T, T+1[$ is almost surely non-empty, given that $Y_T = 1$, but $\mathcal{U} \cap ]G_\infty, G_\infty + 1[ = \emptyset$ by definition. □

**Proposition 6.3.** *If $S$ and $T$ are random variables such that $0 \leq S \leq T \leq \infty$ and $Y$ is continuous on $[S, T[$ (both almost surely) then*

$$Y_t = -W_\bullet(\exp(-1 - t + G_S)) \quad \forall t \in [S, T[$$

*almost surely, where $\bullet \equiv 0$ or $\bullet \equiv -1$ on $[S, T[$.*

**Proof.** Working pathwise, assume $S < T$ and note that, almost surely for all $n \geq 1$, there exists $T_n \in [S, S + 1/n]$ such that $Y_{T_n} \neq 1$ (otherwise $Y \equiv 1$ on $[S, S + 1/n]$, contradicting the fact that $\mathcal{U}$ almost surely has zero Lebesgue measure). Let

$$A := \{R \in ]T_n, T] \colon Y \neq 1 \text{ on } [T_n, R[\};$$

since $Y_{T_n} \neq 1$, the right-continuity of $Y$ at $T_n$ implies that $A$ is non-empty. Furthermore, $R_\infty := \sup A \in A$: there exists $(R_n)_{n \geq 1} \subseteq A$ such that $R_n \nearrow R_\infty$ and $Y \neq 1$ on $\bigcup_{n \geq 1} [T_n, R_n[ = [T_n, R_\infty[$.

If $R \in A$ then, working as in [7], Proof of Theorem 24, it follows that $Y$ is continuously differentiable on $[T_n, R[$ (taking the right derivative at $T_n$) with $Y' = Y/(Y-1)$ there. Hence, by [7], Lemma 25,

$$Y_t = -W_\bullet(-Y_{T_n} \exp(-t + T_n - Y_{T_n})) = -W_\bullet(-\exp(-1 - t + G_{T_n}))$$

for all $t \in [T_n, R[$, where $\bullet \equiv -1$ or $\bullet \equiv 0$. In particular, $Y_{R-} \neq 1$, so if $R_\infty < T$ then $Y$ is continuous at $R_\infty$ and $Y_{R_\infty} \neq 1$, but then there exists $\Delta > 0$ such that $R_\infty + \Delta < T$ and $Y \neq 1$ on $[R_\infty, R_\infty + \Delta[$, contradicting the definition of $R_\infty$. Thus $Y$ has the desired form on $[T_n, T[$; letting $n \to \infty$, so that $T_n \searrow S$, gives the result. $\square$

**Corollary 6.4.** *If $T$ is a random variable such that $Y_T = 1$ almost surely then there exists a sequence $(T_n)_{n \geq 1}$ of random variables such that $T_n \nearrow T$ and $\Delta Y_{T_n} \neq 0$ almost surely.*

**Proof.** Let $T_n := \sup\{t \in \,]0, T]\colon |\Delta Y_t| > 1/(n+1)\}$ for all $n \geq 1$; the sequence $(T_n)_{n \geq 1}$ is increasing, with each $T_n$ almost surely finite and such that $\Delta Y_{T_n} \neq 0$ (since $Y$ has càdlàg paths, so only finitely many jumps of magnitude strictly greater than $1/(n+1)$ on any bounded interval). If $S := \lim_{n \to \infty} T_n$ then $Y$ is continuous on $[S, T[$ and Proposition 6.3 implies that $S = T$ almost surely, as required. $\square$

## 7. Local time

This section is heavily influenced by [17], Section IV.6, hence the proofs are only sketched. Thanks to Theorem 5.5, the results may also be deduced simply from the corresponding properties of Azéma's martingale (except, perhaps, for (21)).

**Definition 7.1.** *Let $\mathcal{P}$ denote the predictable $\sigma$-algebra on $\mathbb{R}_+ \times \Omega$. Recall (see [22], Section I.6, for example) that there exists a $\mathcal{B}(\mathbb{R}) \otimes \mathcal{P}$-measurable function*

$$L \colon \mathbb{R} \times \mathbb{R}_+ \times \Omega \to \mathbb{R}; \qquad (v, t, \omega) \mapsto L_t^v(\omega)$$

*such that, for all $v \in \mathbb{R}$, $L^v$ is a continuous, increasing process with $L_0^v = 0$ and*

$$|Y_t - v| = |v| + \int_0^t \mathrm{sgn}(Y_{s-} - v) \,\mathrm{d}Y_s$$
$$+ \sum_{0 < s \leq t} (|Y_s - v| - |Y_{s-} - v| - \mathrm{sgn}(Y_{s-} - v)\Delta Y_s) + L_t^v \tag{20}$$

*for all $t \geq 0$ almost surely, where $\mathrm{sgn}(x) := \mathbb{1}_{x>0} - \mathbb{1}_{x \leq 0}$ for all $x \in \mathbb{R}$.*

**Remark 7.2.** *Since $X$ is purely discontinuous ([7], Lemma 23), $[Y]^c = [X]^c = 0$; by the occupation-density formula ([17], Corollary 2 to Theorem IV.51), there exists a null set $N \subseteq \Omega$ such that*

$$0 = \int_0^\infty [Y]_t^c(\omega) \,\mathrm{d}t = \int_{-\infty}^\infty \int_0^\infty L_t^v(\omega) \,\mathrm{d}t \,\mathrm{d}v \quad \forall \omega \in \Omega \setminus N,$$

*and so, almost surely, $L^v \equiv 0$ on $\mathbb{R}_+$ for almost all $v \in \mathbb{R}$. The following theorem gives a more exact result.*

**Theorem 7.3.** *If $v \neq 1$ then the local time $L^v = 0$, whereas*

$$\mathbb{E}[L_t^1] = 2 \int_0^t g_\infty(x) \,\mathrm{d}x > 0 \tag{21}$$

*and the random variable $L_t^1$ is not almost surely zero for all $t > 0$.*

**Proof.** If $v = 0$ then (20) implies that

$$|Y_{t+S_0}| = -\int_0^{S_0} dY_s + \int_{S_0}^{t+S_0} dY_s + 2 \sum_{0 < s \le t+S_0} \mathbb{1}_{Y_{s-}=0} \Delta Y_s + L^0_{t+S_0}$$
$$= -1 + Y_{t+S_0} - 1 + 2 + L^0_{t+S_0}$$

for all $t \ge 0$, so $L^0 = 0$. (The first equality uses the local character of the stochastic integral ([17], Corollary to Theorem II.18).) If $v \notin \{0, 1\}$ then the set $\{s > 0: Y_{s-} = Y_s = v\}$ is countable and the claim follows as it does in [17], Proof of Theorem IV.63. For the remaining case, observe that the Meyer–Tanaka–Itô formula (or just [17], Theorem IV.49) yields, for all $t \ge 0$, the identity

$$(Y_t - 1)^+ = \int_0^t \mathbb{1}_{Y_{s-}>1} dY_s + \frac{1}{2} L^1_t.$$

Since

$$\mathbb{E}\left[\int_0^t \mathbb{1}_{Y_{s-}>1} ds\right] = \mathbb{E}\left[\int_0^t \mathbb{1}_{Y_s>1} ds\right] = \int_0^t \mathbb{P}(Y_s > 1) ds,$$

as $\{s > 0: Y_{s-} \ne Y_s\}$ is countable and thus has zero Lebesgue measure, it follows that

$$\mathbb{E}[L^1_t] = 2\mathbb{E}[(Y_t - 1)^+] - 2\int_0^t \mathbb{P}(Y_s > 1) ds.$$

For all $t \ge 0$ and $x \ge 0$, let $F_{Y_t}(x) := \mathbb{P}(Y_t \le x)$; Lemma 2.3 implies that

$$\mathbb{E}[(Y_t - 1)^+] = \int_1^\infty (x-1) dF_{Y_t}(x) = \frac{1}{\pi} \int_1^{b(t)} \operatorname{Im} \frac{x-1}{W_{-1}(-xe^{t-x})} dx = \int_0^t b(t-y) g_\infty(y) dy,$$

using the substitution $x = b(t - y)$, and similarly

$$\int_0^t \mathbb{P}(Y_s > 1) ds = \frac{1}{\pi} \int_0^t \int_1^{b(s)} \operatorname{Im} \frac{1}{W_{-1}(-xe^{s-x})} dx \, ds$$
$$= \int_0^t \int_0^s b'(s-y) g_\infty(y) dy \, ds = \int_0^t (b(t-y) - 1) g_\infty(y) dy.$$

Combining these calculations yields (21). □

**Definition 7.4.** *A semimartingale $R$ has* locally summable jumps *(or satisfies Hypothesis A, in the terminology of [17]) if*

$$\sum_{0 < s \le t} |\Delta R_s| < \infty \quad \text{almost surely } \forall t > 0.$$

**Corollary 7.5.** *The martingale $X$ does not have locally summable jumps.*

**Proof.** Suppose for contradiction that $X$ (and so $Y$) has locally summable jumps. By [17], Theorem IV.56, there exists a $\mathcal{B}(\mathbb{R}) \otimes \mathcal{P}$-measurable function

$$\widetilde{L}: \mathbb{R} \times \mathbb{R}_+ \times \Omega \to \mathbb{R}_+; \qquad (v, t, \omega) \mapsto \widetilde{L}^v_t(\omega)$$

such that $(v, t) \mapsto \widetilde{L}^v_t(\omega)$ is jointly right continuous in $v$ and continuous in $t$ for all $\omega \in \Omega$ and, for all $v \in \mathbb{R}$, $\widetilde{L}^v = L^v$. This is, however, readily seen to contradict Theorem 7.3. □

**Acknowledgments**

This work was begun while the author was a member of the Institut Camille Jordan, Université Claude Bernard Lyon 1, thanks to financial support provided by the European Community's Human Potential Programme under contract HPRN-CT-2002-00279, QP-Applications. It was completed at UCC where the author is an Embark Postdoctoral Fellow funded by the Irish Research Council for Science, Engineering and Technology. Professors Stéphane Attal and Martin Lindsay asked questions which the above attempts to answer; Professor Michel Émery provided several extremely helpful comments on previous drafts. Thanks are extended to the anonymous referee and to the associate editor for their thoughtful and constructive criticism of an earlier version of this paper.

**Appendix A. A Poisson limit theorem**

The following theorem must be well known, but a reference for it (or a version with weaker hypotheses) has proved elusive.

**Theorem A.1.** *For all $n \geq 1$ let $(x_{n,m})_{m=1}^n$ be a collection of independent, identically distributed random variables. If there exists $\lambda > 0$ such that*

$$\lim_{n\to\infty} n\mathbb{E}[x_{n,1}^k] = \lambda \quad \forall k \in \mathbb{N},$$

*then $x_{n,1} + \cdots + x_{n,n}$ converges in distribution to a Poisson law with mean $\lambda$.*

**Proof.** If $n \geq 1$ and $\theta \in \mathbb{R}$ then

$$\left|\mathbb{E}[\exp(\mathrm{i}\theta(x_{n,1} + \cdots + x_{n,n}))] - \left(1 + \frac{\lambda}{n}(\mathrm{e}^{\mathrm{i}\theta} - 1)\right)^n\right| \leq n\left|\mathbb{E}[\mathrm{e}^{\mathrm{i}\theta x_{n,1}}] - 1 - \frac{\lambda(\mathrm{e}^{\mathrm{i}\theta} - 1)}{n}\right|\left(1 + \frac{2\lambda}{n}\right)^{n-1},$$

using the fact that $|z^n - w^n| \leq n|z-w|\max_{1 \leq k \leq n}\{|z|^{k-1}|w|^{n-k}\}$ for all $z$, $w \in \mathbb{C}$. Furthermore, because $|\mathrm{e}^{\mathrm{i}\theta} - \sum_{k=0}^{2p-1}(\mathrm{i}\theta)^k/k!| \leq \theta^{2p}/(2p)!$ for all $\theta \in \mathbb{R}$ and $p \geq 1$,

$$n\left|\mathbb{E}[\mathrm{e}^{\mathrm{i}\theta x_{n,1}}] - 1 - \frac{\lambda(\mathrm{e}^{\mathrm{i}\theta} - 1)}{n}\right|$$

$$\leq n\left|\mathbb{E}\left[\mathrm{e}^{\mathrm{i}\theta x_{n,1}} - \sum_{k=0}^{2p-1}\frac{(\mathrm{i}\theta x_{n,1})^k}{k!}\right]\right| + \sum_{k=1}^{2p-1}\frac{|\theta|^k}{k!}|n\mathbb{E}[x_{n,1}^k] - \lambda| + \lambda\left|\mathrm{e}^{\mathrm{i}\theta} - \sum_{k=0}^{2p-1}\frac{(\mathrm{i}\theta)^k}{k!}\right|$$

$$\leq \frac{|\theta|^{2p}(n\mathbb{E}[x_{n,1}^{2p}] + \lambda)}{(2p)!} + \sum_{k=1}^{2p-1}\frac{|\theta|^k}{k!}|n\mathbb{E}[x_{n,1}^k] - \lambda|.$$

Since $(1 + 2\lambda/n)^{n-1} \to \mathrm{e}^{2\lambda}$ as $n \to \infty$, this sequence is bounded by some constant $C$. Fix $\varepsilon > 0$, choose $p \geq 1$ such that $2|\theta|^{2p}\lambda/(2p)! < \varepsilon/(2C)$ and choose $n_0$ such that

$$\frac{|\theta|^k}{k!}|n\mathbb{E}[x_{n,1}^k] - \lambda| < \frac{\varepsilon}{4pC} \quad \forall n \geq n_0, k = 1, \ldots, 2p;$$

the previous working shows that

$$\left|\mathbb{E}[\exp(\mathrm{i}\theta(x_{n,1} + \cdots + x_{n,n}))] - \left(1 + \frac{\lambda}{n}(\mathrm{e}^{\mathrm{i}\theta} - 1)\right)^n\right| < \frac{2|\theta|^{2p}\lambda C}{(2p)!} + \frac{\varepsilon}{4p} + (2p-1)\frac{\varepsilon}{4p} < \varepsilon \quad \forall n \geq n_0.$$

Hence

$$\lim_{n\to\infty} \mathbb{E}[\exp(i\theta(x_{n,1}+\cdots+x_{n,n}))] = \lim_{n\to\infty}\left(1+\frac{\lambda}{n}(e^{i\theta}-1)\right)^n = \exp(\lambda(e^{i\theta}-1)),$$

and the result follows from the continuity theorem for characteristic functions ([8], Theorem 26.3). □

**Remark A.2.** It follows from the working above that, if $m \geq 1$ and $\theta \in \mathbb{R}$,

$$\mathbb{E}[e^{i\theta x_{n,m}}] = 1 + \left(\frac{\lambda}{n}\right)(e^{i\theta}-1) + o\left(\frac{1}{n}\right) = \mathbb{E}[e^{i\theta b_n}] + o\left(\frac{1}{n}\right) \to 1$$

as $n \to \infty$, where $\mathbb{P}(b_n = 0) = 1 - \lambda/n$ and $\mathbb{P}(b_n = 1) = \lambda/n$. Thus $x_{n,m}$ converges to 0 in distribution, and so in probability, as $n \to \infty$, which explains why this result is a "law of small numbers".

## Appendix B. The probability density function $g_\infty$

**Proposition B.1.** *The function*

$$g_\infty : \mathbb{R} \to \mathbb{R}_+; \qquad x \mapsto \mathbb{1}_{x \geq 0} \frac{1}{\pi} \operatorname{Im} \frac{1}{W_{-1}(-e^{-1+x})}$$

*has a global maximum* $g_\infty(x_0) \approx 0.2306509575$ *at* $x_0 \approx 0.7376612533$, *is strictly increasing on* $[0, x_0]$, *is strictly decreasing on* $[x_0, \infty[$ *with* $\lim_{x\to\infty} g_\infty(x) = 0$,

$$\int_0^\infty g_\infty(x)\,\mathrm{d}x = 1 \quad \text{and} \quad \int_0^\infty x g_\infty(x)\,\mathrm{d}x = \infty.$$

**Proof.** Let $W_{-1}(-e^{-1+x}) = u(x) + iv(x)$ for all $x \geq 0$, where $u(x) \in \mathbb{R}$ and $v(x) \in ]-\pi, 0]$. Since

$$(u+iv)\exp(u+iv) = -\exp(-1+x) \iff \begin{cases} e^u(u\cos v - v\sin v) = -e^{-1+x}, \\ u\sin v + v\cos v = 0, \end{cases}$$

if $v = 0$ then $ue^u = -e^{-1+x}$, which has no solution for $x > 0$, so $v = 0$ if and only if $x = 0$. Suppose henceforth that $x > 0$; note that $u = -v\cot v$,

$$e^{-v\cot v}(-v\cos v \cot v - v\sin v) = -e^{-1+x} \iff x = 1 - v\cot v + \log(v\operatorname{cosec} v)$$

and $\pi g_\infty(x) = -v/(u^2+v^2) > 0$. Observe that

$$\frac{\mathrm{d}u}{\mathrm{d}v} = -\cot v + v\operatorname{cosec}^2 v = \frac{1}{\sin^2 v}(v - \sin v \cos v) < \frac{1}{\sin^2 v}(0 - \sin 0 \cos 0) = 0,$$

because $(\mathrm{d}/\mathrm{d}v)(v - \sin v \cos v) = 1 - \cos 2v > 0$, and

$$\frac{\mathrm{d}x}{\mathrm{d}v} = \frac{1}{v} - 2\cot v + v\operatorname{cosec}^2 v = \frac{1}{v}(1-v\operatorname{cosec} v)^2 - \frac{2}{\sin v}(\cos v - 1) < 0,$$

so $u$ is a strictly increasing function of $x$. As $u(0) = -1$, $u$ takes its values in $[-1, \infty[$; as $v(0) = 0$, letting $x = 1 - v\cot v + \log(v\operatorname{cosec} v) \to \infty$ shows that $v \to -\pi$ (since this function of $v$ is bounded on any proper subinterval of $]-\pi, 0[$) and therefore $u \to \infty$ as $x \to \infty$. (In particular, $|g_\infty(x)| \leq 1/u^2 \to 0$ as $x \to \infty$.) Since $u$ is continuous, strictly increasing and maps $[0, \infty[$ to $[-1, \infty[$, there exists $x_0$ such that $u(x_0) = -1/2$. Moreover,

$$g'_\infty(x) = \operatorname{Im} \frac{\mathrm{d}}{\mathrm{d}x} \frac{1}{W_{-1}(-e^{-1+x})} = \operatorname{Im} \frac{-1}{(u+iv)(1+u+iv)} = \frac{v(2u+1)}{(u^2+v^2)((1+u)^2+v^2)},$$

so $g'_\infty > 0$ on $]0, x_0[$ and $g'_\infty < 0$ on $]x_0, \infty[$. (The approximate values for $x_0$ and $g_\infty(x_0)$ were determined with the use of Maple.)

For the integrals, the substitution $x = v$ gives that

$$\pi \int_0^\infty g_\infty(x)\,dx = \int_{-\pi}^0 \frac{\sin^2 v}{v}\left(\frac{1}{v} - 2\cot v + v\csc^2 v\right) dv$$

$$= \pi + \int_{-\pi}^0 \left(\frac{\sin^2 v}{v^2} - \frac{\sin 2v}{v}\right) dv = \pi + \left[-\frac{\sin^2 v}{v}\right]_{-\pi}^0 = \pi,$$

as required. Finally, if $\varepsilon \in ]0, \pi/2[$,

$$\pi \int_0^\infty x g_\infty(x)\,dx = \int_{-\pi}^0 (1 - v\cot v + \log(v\csc v))\left(\frac{\sin^2 v}{v^2} - \frac{\sin 2v}{v} + 1\right) dv$$

$$\geq \int_{-\pi+\varepsilon}^{-\pi/2} -v\cot v\,dv \geq \frac{\pi}{2} \int_\varepsilon^{\pi/2} \cot w\,dw = -\log\sin\varepsilon \to \infty$$

as $\varepsilon \to 0+$. $\qquad\square$

**Remark B.2.** *It follows from Propositions B.1 and 3.1 that the distribution of $G_\infty$ is unimodal with mode $x_0$, that is, $t \mapsto \mathbb{P}(G_\infty \leq t)$ is convex on $]-\infty, x_0[$ and concave on $]x_0, \infty[$.*

## Appendix C. An auxiliary calculation

**Lemma C.1.** *If $f_J$ is as defined in Proposition 4.6 then*

$$\pi f_J(t) = \operatorname{Im} \frac{1}{1 + W_{-1}(-e^{-1+t})} \sim \frac{1}{\sqrt{2t}} \quad \text{as } t \to 0+$$

*and $f_J$ is strictly decreasing on $]0, \infty[$.*

**Proof.** For all $t \geq 0$, let $p := -\sqrt{2(1-e^t)} = -i\sqrt{2t} + O(t^{3/2})$ as $t \to 0+$; recall that

$$-W_{-1}(-e^{-1+t}) = 1 - p + O(p^2) = 1 + i\sqrt{2t} + O(t)$$

as $t \to 0+$, by [10], (4.22), and this gives the first result. For the next claim, if $t > 0$ and $W_{-1}(-e^{-1+t}) = -v\cot v + iv$, where $v \in ]-\pi, 0[$, then

$$\pi f'_J(t) = \operatorname{Im} \frac{-W_{-1}(-e^{-1+t})}{(1 + W_{-1}(-e^{-1+t}))^3} = \frac{((3 - 2v\cot v)v^2 \csc^2 v - 1)v}{((1 - v\cot v)^2 + v^2)^3}.$$

The result follows if

$$(3 - 2v\cot v)v^2\csc^2 v - 1 > 0 \iff (v^2 - \sin^2 v)\sin v + 2v^2(\sin v - v\cos v) < 0$$

for all $v \in ]-\pi, 0[$, but since $\sin^2 v < v^2$ and $\sin v < v\cos v$ for such $v$, this is clear. $\qquad\square$

**Proposition C.2.** *If $D := \{(t, x) \in \mathbb{R}_+^2 : a(t) \leq x \leq b(t)\}$,*

$$f: D \to \mathbb{R}_+; \qquad (t, x) \mapsto \operatorname{Im} \frac{1}{W_{-1}(-xe^{t-x})},$$

$$F: D \to \mathbb{R}_+; \qquad (t, x) \mapsto \int_{a(t)}^x f(t, y)\,dy$$

*and* $(s, y) \in D^\circ := \{(t, x) \in \mathbb{R}_+^2 : t > 0, a(t) < x < b(t)\}$ *then*

$$F(s, y) + \frac{\partial F}{\partial t}(s, y) = \int_{a(s)}^{y} \operatorname{Im} \frac{1}{1 + W_{-1}(-z e^{s-z})} \, dz. \tag{22}$$

**Proof.** Note first that, since $f$ is continuous, $F$ is well defined. If $h > 0$ then

$$\frac{F(s+h, y) - F(s, y)}{h} = \frac{1}{h} \int_{a(s+h)}^{a(s)} f(s+h, z) \, dz + \int_{a(s)}^{y} \frac{f(s+h, z) - f(s, z)}{h} \, dz$$

and the intermediate-value theorem gives $\zeta_h \in [a(s+h), a(s)]$ such that

$$\frac{1}{h} \int_{a(s+h)}^{a(s)} f(s+h, z) \, dz = \frac{a(s) - a(s+h)}{h} f(s+h, \zeta_h) \to -a'(s) f(s, a(s)) = 0$$

as $h \to 0+$. For all $z \in [a(s), b(s)]$ there exists $\theta_{h,z} \in ]0, 1[$ such that

$$\frac{f(s+h, z) - f(s, z)}{h} = \frac{\partial f}{\partial t}(s + \theta_{h,z} h, z)$$

by the mean-value theorem, since $t \mapsto f(t, z)$ is continuous on $[s, s+h]$ and differentiable on $]s, s+h[$. Let

$$g : D^\circ \to \mathbb{R}_+; \qquad (t, x) \mapsto \frac{\partial f}{\partial t}(t, x) + f(t, x) = \begin{cases} \pi f_J(t - a^{-1}(x)) & \text{if } x \in ]a(t), 1], \\ \pi f_J(t - b^{-1}(x)) & \text{if } x \in [1, b(t)[, \end{cases}$$

where $f_J$ is defined in Proposition 4.6. The continuity of $f$ on $[s, s+1] \times [a(s), y]$ and the dominated-convergence theorem imply that

$$F(s, y) = \int_{a(s)}^{y} f(s, z) \, dz = \lim_{h \to 0+} \int_{a(s)}^{y} f(s + \theta_{h,z} h, z) \, dz,$$

so the right-hand limit in (22) has the correct value if $\int_{a(s)}^{y} g(s, z) \, dz$ exists and

$$\lim_{h \to 0+} \int_{a(s)}^{y} g(s + \theta_{h,z} h, z) \, dz = \int_{a(s)}^{y} g(s, z) \, dz.$$

Fix $r \in ]0, s[$ such that $y > a(r)$ and note that $g$ is continuous on $[s, s+1] \times [a(r), y]$, so the dominated-convergence theorem implies that

$$\lim_{h \to 0+} \int_{a(r)}^{y} g(s + \theta_{h,z} h, z) \, dz = \int_{a(r)}^{y} g(s, z) \, dz.$$

Next, note that if $z \in ]a(s), a(r)]$ and $h \to 0+$ then

$$g(s + \theta_{h,z} h, z) = \pi f_J(s + \theta_{h,z} h - a^{-1}(z)) \nearrow \pi f_J(s - a^{-1}(z)) = g(s, z),$$

because $f_J$ is strictly decreasing, by Lemma C.1. The first half of the result now follows from the monotone-convergence theorem, once it is known that $\int_{a(s)}^{a(r)} g(s, z) \, dz$ exists. However,

$$\int_{a(s)}^{a(r)} g(s, z) \, dz = \pi \int_{s}^{r} f_J(s - u) a'(u) \, du = -\pi \int_{0}^{s-r} f_J(t) a'(s - t) \, dt < \infty,$$

since, by Lemma C.1, $\pi f_J(t) \sim 1/\sqrt{2t}$ as $t \to 0+$, $f_J$ is continuous on $]0, s - r]$ and $a'$ is continuous on $[r, s]$.

Now suppose that $h < 0$ is such that $s + h > 0$ and $b(s+h) > y > a(s+h)$. Then

$$\frac{F(s+h,y) - F(s,y)}{h} = \int_{a(s+h)}^{y} \frac{f(s+h,z) - f(s,z)}{h}\,dz - \frac{1}{h}\int_{a(s)}^{a(s+h)} f(s,z)\,dz$$

and the second term tends to 0 as $h \to 0-$. If $z \in [a(s+h), b(s+h)]$ then $t \mapsto f(t,z)$ is continuous on $[s+h, s]$ and differentiable on $]s+h, s[$, so there exists $\theta_{h,z} \in\,]0,1[$ such that

$$\frac{f(s+h,z) - f(s,z)}{h} = \frac{\partial f}{\partial t}(s + \theta_{h,z} h, z).$$

Furthermore, as $f$ is continuous, so bounded, on the compact set $D \cap ([0,s] \times \mathbb{R}_+)$, the dominated-convergence theorem implies that

$$F(s,y) = \lim_{h \to 0-} \int_{a(s+h)}^{y} f(s + \theta_{h,z} h, z)\,dz$$

and the result follows if

$$\lim_{h \to 0-} \int_{a(s+h)}^{y} g(s + \theta_{h,z} h, z)\,dz = \int_{a(s)}^{y} g(s,z)\,dz.$$

Fix $0 < r_1 < r_2 < s$ such that $a(r_1) < y$ and note that $g$ is continuous on $[r_2, s] \times [a(r_1), y]$, so bounded there, and the dominated-convergence theorem implies that

$$\int_{a(r_1)}^{y} g(s + \theta_{h,z} h, z)\,dz \to \int_{a(r_1)}^{y} g(s,z)\,dz$$

as $h \to 0-$. A final application of the monotone-convergence theorem completes the result, since if $h < 0$ is such that $r_2 < s + h$ then, letting $h \to 0-$,

$$\mathbb{1}_{z \in [a(s+h), a(r_1)]} g(s + \theta_{h,z} h, z) = \mathbb{1}_{z \in [a(s+h), a(r_1)]} \pi f_J(s + \theta_{h,z} h - a^{-1}(z))$$
$$\nearrow \mathbb{1}_{z \in ]a(s), a(r_1)]} \pi f_J(s - a^{-1}(z))$$
$$= \mathbb{1}_{z \in ]a(s), a(r_1)]} g(s,z). \qquad \square$$

### Appendix D. A pair of Laplace transforms

**Theorem D.1.** *If $g_\infty$ is as defined in Proposition* 3.1 *and $f_J$ is as defined in Proposition* 4.6 *then their Laplace transforms are as follows:*

$$\widehat{g_\infty}(p) = \frac{e^{-p} p^p}{\Gamma(p+2)} \quad \text{and} \quad \widehat{f_J}(p) = (p+1)\widehat{g_\infty}(p) = \frac{e^{-p} p^p}{\Gamma(p+1)}, \qquad (23)$$

*where* $\Gamma\colon p \mapsto \int_0^\infty z^{p-1} e^{-z}\,dz$ *is the gamma function.*

**Proof.** Let

$$f_1(t) := \frac{1}{\pi}\int_{a(t)}^{b(t)} \operatorname{Im}\frac{1}{W_{-1}(-y e^{t-y})}\,dy \quad \forall t \geq 0.$$

Splitting the interval $[a(t), b(t)]$ at 1 and using the substitutions $y = a(t-x)$ and $y = b(t-x)$, as appropriate,

$$f_1(t) = \frac{1}{\pi}\int_0^t \operatorname{Im}\left(\frac{1}{W_{-1}(-e^{-1+x})}\right) c(t-x)\,dx = (g_\infty \star c)(t),$$

where $\star$ denotes convolution of functions on $\mathbb{R}_+$ and $c$ is as in Definition 2.2. Furthermore,

$$\widehat{c}(p) := \int_0^\infty c(x)\mathrm{e}^{-px}\,\mathrm{d}x = \int_0^\infty b'(x)\mathrm{e}^{-px}\,\mathrm{d}x - \int_0^\infty a'(x)\mathrm{e}^{-px}\,\mathrm{d}x$$

$$= \int_1^\infty \mathrm{e}^{-p(-1+y-\log y)}\,\mathrm{d}y + \int_0^1 \mathrm{e}^{-p(-1+y-\log y)}\,\mathrm{d}y$$

$$= \mathrm{e}^p \int_0^\infty \left(\frac{z}{p}\right)^p \mathrm{e}^{-z} p^{-1}\,\mathrm{d}z.$$

The second line follows from the substitutions $y = b(x)$ and $y = a(x)$. Thus, since $f_1(t) = \mathbb{P}(Y_t > 0) = 1 - \mathrm{e}^{-t}$,

$$\widehat{g_\infty}(p) = \frac{\widehat{f_1}(p)}{\widehat{c}(p)} = \frac{1}{p(p+1)} \frac{\mathrm{e}^{-p}p^{p+1}}{\Gamma(p+1)} = \frac{\mathrm{e}^{-p}p^p}{\Gamma(p+2)},$$

as claimed. If

$$f_2(t) := \frac{1}{\pi} \int_{a(t)}^{b(t)} \operatorname{Im} \frac{1}{1 + W_{-1}(-z\mathrm{e}^{t-z})}\,\mathrm{d}y \quad \forall t > 0,$$

then, working as above, $f_2 = f_J \star c$. Moreover, since $f_2 = f_1 + f_1'$ (by the working in the proof of Theorem 4.5), it follows that $\widehat{f_2}(p) = (p+1)\widehat{f_1}(p)$ and

$$\widehat{f_J}(p) = \frac{\widehat{f_2}(p)}{\widehat{c}(p)} = \frac{(p+1)\widehat{f_1}(p)}{\widehat{c}(p)} = \frac{\mathrm{e}^{-p}p^p}{\Gamma(p+1)}.$$

$\square$

**Remark D.2.** *The substitution $x = 1 - v \cot v + \log(v \operatorname{cosec} v)$ yields the identity*

$$\mathrm{e}^p \widehat{f_J}(p) = \frac{1}{\pi} \int_0^\pi \left(\frac{\sin v}{v}\right)^p \exp(pv \cot v)\,\mathrm{d}v; \tag{24}$$

*it should be possible to verify directly that the right-hand side of (24) equals $p^p/\Gamma(p+1)$. (This would give independent proof that*

$$A \mapsto \mathbb{1}_{0 \in A} + \frac{1}{\pi} \int_{A \cap [a(t), b(t)]} \operatorname{Im} \frac{1}{W_{-1}(-y\mathrm{e}^{t-y})}\,\mathrm{d}y$$

*and*

$$A \mapsto \frac{1}{\pi} \int_{A \cap [a(t), b(t)]} \operatorname{Im} \frac{1}{1 + W_{-1}(-z\mathrm{e}^{t-z})}\,\mathrm{d}z$$

*are probability measures on $\mathcal{B}(\mathbb{R})$.)*

### References

[1] S. Attal. The structure of the quantum semimartingale algebras. *J. Operator Theory* **46** (2001) 391–410. MR1870414
[2] S. Attal and A. C. R. Belton. The chaotic-representation property for a class of normal martingales. *Probab. Theory Related Fields* **139** (2007) 543–562. MR2322707
[3] J. Azéma. Sur les fermés aléatoires. *Séminaire de Probabilités XIX* 397–495. J. Azéma and M. Yor (Eds). *Lecture Notes in Math.* **1123**. Springer, Berlin, 1985. MR0889496
[4] J. Azéma and M. Yor. Étude d'une martingale remarquable. *Séminaire de Probabilités XXIII* 88–130. J. Azéma, P.-A. Meyer and M. Yor (Eds). *Lecture Notes in Math.* **1372**. Springer, Berlin, 1989. MR1022900
[5] A. C. R. Belton. An isomorphism of quantum semimartingale algebras. *Q. J. Math.* **55** (2004) 135–165. MR2068315

[6] A. C. R. Belton. A note on vacuum-adapted semimartingales and monotone independence. In *Quantum Probability and Infinite Dimensional Analysis XVIII. From Foundations to Applications*, 105–114. M. Schürmann and U. Franz (Eds), World Scientific, Singapore, 2005. MR2211883

[7] A. C. R. Belton. The monotone Poisson process. In *Quantum Probability* 99–115. M. Bożejko, W. Młotkowski and J. Wysoczański (Eds). Banach Center Publications **73**, Polish Academy of Sciences, Warsaw, 2006.

[8] P. Billingsley. *Probability and Measure*, 3rd edition. Wiley, New York, 1995. MR1324786

[9] C. S. Chou. Caractérisation d'une classe de semimartingales. *Séminaire de Probabilités XIII* 250–252. C. Dellacherie, P.-A. Meyer and M. Weil (Eds). *Lecture Notes in Math.* **721**. Springer, Berlin, 1979. MR0544798

[10] R. M. Corless, G. H. Gonnet, D. E. G. Hare, D. J. Jeffrey and D. E. Knuth. On the Lambert $W$ function. *Adv. Comput. Math.* **5** (1996) 329–359. MR1414285

[11] F. Delbaen and W. Schachermayer. The fundamental theorem of asset pricing for unbounded stochastic processes. *Math. Ann.* **312** (1998) 215–250. MR1671792

[12] M. Émery. Compensation de processus à variation finie non localement intégrables. *Séminaire de Probabilités XIV* 152–160. J. Azéma and M. Yor (Eds). *Lecture Notes in Math.* **784**. Springer, Berlin, 1980. MR0580120

[13] M. Émery. On the Azéma martingales. *Séminaire de Probabilités XXIII* 66–87. J. Azéma, P.-A. Meyer and M. Yor (Eds). *Lecture Notes in Math.* **1372**. Springer, Berlin, 1989. MR1022899

[14] M. Émery. Personal communication, 2006.

[15] R. L. Graham, D. E. Knuth and O. Patashnik. *Concrete Mathematics*, 2nd edition. Addison-Wesley, Reading, MA, 1994. MR1397498

[16] N. Muraki. Monotonic independence, monotonic central limit theorem and monotonic law of small numbers. *Infin. Dimens. Anal. Quantum Probab. Relat. Top.* **4** (2001) 39–58. MR1824472

[17] P. Protter. *Stochastic Integration and Differential Equations. A New Approach*. Springer, Berlin, 1990. MR1037262

[18] L. C. G. Rogers and D. Williams. *Diffusions, Markov Processes and Martingales. Volume 1: Foundations*, 2nd edition. Cambridge University Press, Cambridge, 2000. MR1796539

[19] W. Rudin. *Real and Complex Analysis*, 3rd edition. McGraw-Hill, New York, 1987. MR0924157

[20] R. Speicher. A new example of "independence" and "white noise". *Probab. Theory Related Fields* **84** (1990) 141–159. MR1030725

[21] C. Stricker. Représentation prévisible et changement de temps. *Ann. Probab.* **14** (1986) 1070–1074. MR0841606

[22] C. Stricker and M. Yor. Calcul stochastique dépendant d'un paramètre. *Z. Wahrsch. Verw. Gebiete* **45** (1978) 109–133. MR0510530

[23] G. Taviot. Martingales et équations de structure: étude géométrique. Thèse, Université Louis Pasteur Strasbourg 1, 1999. MR1736397

[24] S. J. Taylor. The $\alpha$-dimensional measure of the graph and set of zeros of a Brownian path. *Proc. Cambridge Philos. Soc.* **51** (1955) 265–274. MR0074494